\documentclass[12pt]{article}
\usepackage{latexsym,amsfonts,amsmath,amssymb}
\newtheorem{theorem}{Theorem}
\newtheorem{corollary}[theorem]{Corollary}
\newtheorem{sublemma}{Lemma}[theorem]
\newtheorem{lemma}[theorem]{Lemma}
\newtheorem{question}[theorem]{Question}
\newtheorem{observation}[theorem]{Observation}
\newtheorem{claim}[theorem]{Claim}
\newtheorem{conjecture}[theorem]{Conjecture}
\newtheorem{definition}[theorem]{Definition}
\newtheorem{remark}[theorem]{Remark}
\newtheorem{example}[theorem]{Example}
\def\Theorem #1.#2 #3\par{\setbox1=\hbox{#1}\ifdim\wd1=0pt
   \begin{theorem}{\rm #2} #3\end{theorem}\else
   \newtheorem{#1}[theorem]{#1}\begin{#1}\label{#1}{\rm #2} #3\end{#1}\fi}
\def\Corollary #1.#2 #3\par{\setbox1=\hbox{#1}\ifdim\wd1=0pt
   \begin{corollary}{\rm #2} #3\end{corollary}\else
   \newtheorem{#1}[theorem]{#1}\begin{#1}\label{#1}{\rm #2} #3\end{#1}\fi}
\def\Lemma #1.#2 #3\par{\setbox1=\hbox{#1}\ifdim\wd1=0pt
   \begin{lemma}{\rm #2} #3\end{lemma}\else
   \newtheorem{#1}[theorem]{#1}\begin{#1}\label{#1}{\rm #2} #3\end{#1}\fi}
\def\SubLemma #1.#2 #3\par{\setbox1=\hbox{#1}\ifdim\wd1=0pt
   \begin{sublemma}{\rm #2} #3\end{sublemma}\else
   \newtheorem{#1}{#1}[theorem]\begin{#1}\label{#1}{\rm #2} #3\end{#1}\fi}
\def\Question #1.#2 #3\par{\setbox1=\hbox{#1}\ifdim\wd1=0pt
   \begin{question}{\rm #2} #3\end{question}\else
   \newtheorem{#1}[theorem]{#1}\begin{#1}\label{#1}{\rm #2} #3\end{#1}\fi}
\def\Observation #1.#2 #3\par{\setbox1=\hbox{#1}\ifdim\wd1=0pt
   \begin{observation}{\rm #2} #3\end{observation}\else
   \newtheorem{#1}[theorem]{#1}\begin{#1}\label{#1}{\rm #2} #3\end{#1}\fi}
\def\Claim #1.#2 #3\par{\setbox1=\hbox{#1}\ifdim\wd1=0pt
   \begin{claim}{\rm #2} #3\end{claim}\else
   \newtheorem{#1}[theorem]{#1}\begin{#1}\label{#1}{\rm #2} #3\end{#1}\fi}
\def\Conjecture #1.#2 #3\par{\setbox1=\hbox{#1}\ifdim\wd1=0pt
   \begin{conjecture}{\rm #2} #3\end{conjecture}\else
   \newtheorem{#1}[theorem]{#1}\begin{#1}\label{#1}{\rm #2} #3\end{#1}\fi}
\def\Definition #1.#2 #3\par{\setbox1=\hbox{#1}\ifdim\wd1=0pt
   \begin{definition}{\rm #2} #3\end{definition}\else
   \newtheorem{#1}[theorem]{#1}\begin{#1}\label{#1}{\rm #2} #3\end{#1}\fi}
\def\Remark #1.#2 #3\par{\setbox1=\hbox{#1}\ifdim\wd1=0pt
   \begin{remark}{\rm #2} #3\end{remark}\else
   \newtheorem{#1}[theorem]{#1}\begin{#1}\label{#1}{\rm #2} #3\end{#1}\fi}
\def\Example #1.#2 #3\par{\setbox1=\hbox{#1}\ifdim\wd1=0pt
   \begin{example}{\rm #2} #3\end{example}\else
   \newtheorem{#1}[theorem]{#1}\begin{#1}\label{#1}{\rm #2} #3\end{#1}\fi}
\def\QuietTheorem #1.#2 #3\par{\setbox1=\hbox{#1}\medskip\ifdim\wd1=0pt
   \proclaim{Theorem {\rm #2}}{#3}\else\proclaim{#1 {\rm #2}}{#3}\fi}
\newcommand{\proclaim}[2]{\medskip\noindent{\bf #1} {\sl#2}\par\medskip\noindent}
\def\Proclaim #1.#2 #3\par{\proclaim{#1 {\rm #2}}{#3}}
\newenvironment{proof}{\noindent}{\kern2pt\QEDbox\par\bigskip}
\def\Proof#1: {\setbox1=\hbox{#1}\ifdim\wd1=0pt\begin{proof}{\bf Proof: }\else\medskip\begin{proof}{\bf #1: }\fi}
\newcommand{\QED}{\end{proof}}
\def\Abstract #1\par{\medskip\begin{quotation}{\singlespaced\footnotesize{\noindent{\bf Abstract.~}#1}}\end{quotation}\medskip}
\def\Title #1\par{\title{#1}\maketitle}
\def\Author #1\par{\author{#1}}
\def\Acknowledgement#1\par{\thanks{#1}}
\def\Chapter #1\par{\chapter{#1}}
\def\Section #1\par{\section{#1}}
\def\SubSection #1\par{\subsection{#1}}
\def\SubSubSection #1\par{\subsubsection{#1}}
\def\MidTitle #1\par{\bigskip\goodbreak\centerline{\small\bf #1}\bigskip}

\newcommand{\singlespaced}{\baselineskip=15pt}
\renewcommand{\P}{{\mathbb P}}
\newcommand{\Q}{{\mathbb Q}}
\newcommand{\R}{{\mathbb R}}

\newcommand{\Gtail}{G_{\!\scriptscriptstyle\rm tail}}

\newcommand{\Ptail}{\P_{\!\scriptscriptstyle\rm tail}}

\newcommand{\Pterm}{\P_{\!\scriptscriptstyle\rm term}}
\newcommand{\Gterm}{G_{\!\scriptscriptstyle\rm term}}

\newcommand{\one}{\mathop{1\hskip-3pt {\rm l}}}
\newfont{\msam}{msam10 at 12pt}

\newcommand{\of}{\subseteq}

\newcommand{\set}[1]{\left\{\,{#1}\,\right\}}

\newcommand{\compose}{\circ}
\newcommand{\elesub}{\prec}

\newcommand{\dom}{\mathop{\rm dom}}
\newcommand{\ran}{\mathop{\rm ran}}
\newcommand{\add}{\mathop{\rm add}}
\newcommand{\coll}{\mathop{\rm coll}}

\newcommand{\Cof}{\mathop{\rm Cof}}

\newcommand{\image}{\mathbin{\hbox{\tt\char'42}}}
\newcommand{\plus}{{+}}
\newcommand{\plusplus}{{{+}{+}}}
\newcommand{\plusplusplus}{{{+}{+}{+}}}

\newcommand{\forces}{\Vdash}

\newcommand{\cross}{\times}
\newcommand{\union}{\cup}

\newcommand{\intersect}{\cap}

\newcommand{\smalllt}{\mathrel{\mathchoice{\raise2pt\hbox{$\scriptstyle<$}}{\raise1pt\hbox{$\scriptstyle<$}}{\scriptscriptstyle<}{\scriptscriptstyle<}}}
\newcommand{\smallleq}{\mathrel{\mathchoice{\raise2pt\hbox{$\scriptstyle\leq$}}{\raise1pt\hbox{$\scriptstyle\leq$}}{\scriptscriptstyle\leq}{\scriptscriptstyle\leq}}}
\newcommand{\ltkappa}{{{\smalllt}\kappa}}
\newcommand{\leqkappa}{{{\smallleq}\kappa}}

\newcommand{\leqlambda}{{{\smallleq}\lambda}}

\newcommand{\ltgamma}{{{\smalllt}\gamma}}
\newcommand{\leqeta}{{{\smallleq}\eta}}
\newcommand{\lteta}{{{\smalllt}\eta}}
\newcommand{\leqzeta}{{{\smallleq}\zeta}}

\newcommand{\leqtheta}{{{\smallleq}\theta}}
\newcommand{\lttheta}{{{\smalllt}\theta}}

\newcommand{\leqdelta}{{{\smallleq}\delta}}

\newcommand{\card}[1]{{\left|#1\right|}}

\newcommand{\UnderTilde}[1]{{\setbox1=\hbox{$#1$}\baselineskip=0pt\vtop{\hbox{$#1$}\hbox to\wd1{\hfil$\sim$\hfil}}}{}}
\newcommand{\Undertilde}[1]{{\setbox1=\hbox{$#1$}\baselineskip=0pt\vtop{\hbox{$#1$}\hbox to\wd1{\hfil$\scriptstyle\sim$\hfil}}}{}}
\newcommand{\undertilde}[1]{{\setbox1=\hbox{$#1$}\baselineskip=0pt\vtop{\hbox{$#1$}\hbox to\wd1{\hfil$\scriptscriptstyle\sim$\hfil}}}{}}
\newcommand{\UnderdTilde}[1]{{\setbox1=\hbox{$#1$}\baselineskip=0pt\vtop{\hbox{$#1$}\hbox to\wd1{\hfil$\approx$\hfil}}}{}}
\newcommand{\Underdtilde}[1]{{\setbox1=\hbox{$#1$}\baselineskip=0pt\vtop{\hbox{$#1$}\hbox to\wd1{\hfil\scriptsize$\approx$\hfil}}}{}}

\newcommand{\st}{\mid}

\def\<#1>{\langle\,#1\,\rangle}

\newcommand{\QEDbox}{\fbox{}}
\newcommand{\cp}{\mathop{\rm cp}}

\newcommand{\GCH}{\hbox{\sc gch}}

\newcommand{\factordiagramup}[6]{$$\begin{array}{ccc}
#1&\raise3pt\vbox{\hbox to60pt{\hfill$\scriptstyle
#2$\hfill}\vskip-6pt\hbox{$\vector(4,0){60}$}}&#3\\ \vbox
to30pt{}&\raise22pt\vtop{\hbox{$\vector(4,-3){60}$}\vskip-22pt\hbox
to60pt{\hfill$\scriptstyle #4\qquad$\hfill}}
     &\ \ \lower22pt\hbox{$\vector(0,3){45}$}\ {\scriptstyle #5}\\
\vbox to15pt{}&&#6\\
\end{array}$$}
\newcommand{\factordiagram}[6]{$$\begin{array}{ccc}
#1&&\\ \ \ \raise22pt\hbox{$\vector(0,-3){45}$}\ {\scriptstyle #2}
&\raise22pt\hbox{$\vector(2,-1){90}$}\raise5pt\llap{$\scriptstyle#3$\qquad\quad}&\vbox
to25pt{}\\ #4&\raise3pt\vbox{\hbox to90pt{\hfill$\scriptstyle
#5$\hfill}\vskip-6pt\hbox{$\vector(4,0){90}$}}&#6\\
\end{array}$$}
\newcommand{\df}{\it} % use italic for definition terms

\begin{document}

\author{Arthur W. Apter\thanks{The authors' research has been supported in
part by PSC-CUNY grants from the CUNY Research Foundation and, in
the case of the second author,
from the NSF.}\\
{\normalsize\sc Department of Mathematics}\\
{\normalsize\sc Baruch College of CUNY}\\
{\normalsize\sc New York, NY 10010 USA}\\
{\footnotesize http://math.baruch.cuny.edu/$\sim$apter}\\
{\footnotesize awabb@cunyvm.cuny.edu}\\
\\
Joel David Hamkins\footnotemark[\value{footnote}]\\
\normalsize\sc The City University of New
York\footnote{Specifically,
The College of Staten Island and The CUNY Graduate Center.}\\
\normalsize\sc Carnegie Mellon University\\
{\footnotesize http://jdh.hamkins.org}\\
{\footnotesize jdh@hamkins.org}\\}

\Title Indestructibility and the level-by-level agreement between
strong compactness and supercompactness\thanks{AMS Math Subject
Codes: 03E35, 03E55. Keywords: Supercompact cardinal, strongly
compact cardinal, indestructibility, level-by-level agreement}

\Abstract Can a supercompact cardinal $\kappa$ be Laver
indestructible when there is a level-by-level agreement between
strong compactness and supercompactness? In this article, we show
that if there is a sufficiently large cardinal above $\kappa$,
then no, it cannot. Conversely, if one weakens the requirement
either by demanding less indestructibility, such as requiring
only indestructibility by stratified posets, or less
level-by-level agreement, such as requiring it only on measure
one sets, then yes, it can.

Two important but apparently unrelated results occupy the large
cardinal literature. On the one hand, Laver \cite{Laver78}
famously proved that any supercompact cardinal $\kappa$ can be
made indestructible by $\ltkappa$-directed closed forcing. On the
other hand, Apter and Shelah \cite{ApterShelah495} proved that
all supercompact cardinals can be preserved to a forcing extension
where there is a level-by-level agreement between strong
compactness and supercompactness: specifically, except in special
cases known to be impossible, any cardinal $\gamma$ there is
$\eta$-strongly compact if and only if it is
$\eta$-supercompact.\footnote{See the definition in the paragraph
immediately preceding Observation \ref{Stratified}.} Can these
results be combined? Specifically, we ask:

\Question Open Question. Can a supercompact cardinal be
indestructible when there is a level-by-level agreement between
strong compactness and supercompactness?

In this article, we provide a partial answer to this question,
constraining the possibilities from both above and below. But
alas, our results do not settle the matter, so the question
remains open. What we can prove, specifically, is that if there
is a sufficiently large cardinal above the supercompact cardinal,
then the answer to the question is {\it no}. In particular, there
is at most one supercompact cardinal as in the question; more
exactly, if a cardinal is indestructibly supercompact in the
presence of a level-by-level agreement between strong compactness
and supercompactness, then no larger cardinal $\lambda$ is
$2^\lambda$-supercompact. Conversely, if the requirements in the
question are weakened in any of several ways, asking either for
less indestructibility, replacing it with resurrectibility or
with indestructibility by stratified forcing, or for a weaker
form of level-by-level agreement, demanding that it hold only on
measure one sets, then the answer is {\it yes}. These results are
summarized in the Main Theorem stated below.

\QuietTheorem Main Theorem.
\begin{enumerate}
\item There can be at most one supercompact cardinal as in the
question; indeed, if\/ $\kappa$ is indestructibly supercompact and
there is a level-by-level agreement between strong compactness and
supercompactness, then no cardinal $\lambda$ above $\kappa$ is
$2^\lambda$-supercompact. This same conclusion can be made if we
assume only that $\kappa$ is indestructibly strong, or
indestructibly $\Sigma_2$-reflecting.
\item Conversely, relaxing the notion of indestructibility somewhat,
it is relatively consistent to have in the presence of a
level-by-level agreement between strong compactness and
supercompactness a supercompact cardinal $\kappa$ that is
indestructible by any stratified $\ltkappa$-directed closed
forcing and more. It follows that the supercompactness of
$\kappa$ is resurrectible after any $\ltkappa$-directed closed
forcing.
\item Alternatively, by relaxing the degree of level-by-level
agreement required, it is relatively consistent to have a fully
indestructible supercompact cardinal with a level-by-level
agreement almost everywhere between strong compactness and
supercompactness.
\end{enumerate}

The precise details of these three claims---including
definitions, stronger statements of the results and
corollaries---appear respectively in the three sections of this
article.

We will define the most important notions here. We say that a
supercompact cardinal $\kappa$ is {\df indestructible} when it
remains supercompact after any $\ltkappa$-directed closed
forcing. A forcing notion is {\df $\ltkappa$-directed closed}
when any directed subset of it of size less than $\kappa$ has a
lower bound. If additionally any directed subset of it of size
$\kappa$ has a lower bound, then the forcing notion is {\df
$\leqkappa$-directed closed}. A forcing notion is {\df
$\leqkappa$-closed} if any decreasing chain of length less than
or equal to $\kappa$ has a lower bound. A forcing notion is {\df
$\ltkappa$-distributive} if forcing with it adds no new sequences
over the ground model of length less than $\kappa$. If the
forcing notion adds no new $\kappa$-sequences, then it is {\df
$\leqkappa$-distributive}. A forcing notion is {\df
$\leqkappa$-strategically closed} if in the game of length
$\kappa + 1$ in which two players alternately select conditions
from it to construct a descending $\kappa$-sequence, with the
second player playing at limit stages, the second player has a
strategy that allows her always to continue playing. The
supercompactness of $\kappa$ is {\df resurrectible} after any
$\ltkappa$-directed closed forcing if, after any such forcing
$\Q$ there is further $\ltkappa$-distributive forcing $\R$ such
that $\Q*\R$ preserves the supercompactness of $\kappa$. (We will
occasionally abuse notation and write $x$ when we should more
properly write $\dot x$.) We say that there is a {\df
level-by-level agreement} between strong compactness and
supercompactness when for any two regular cardinals
$\gamma\leq\eta$, the cardinal $\gamma$ is $\eta$-strongly
compact if and only if it is $\eta$-supercompact, unless $\gamma$
is a measurable limit of cardinals which are $\eta$-strongly
compact. We say that $\gamma$ is {\df partially supercompact} if
and only if $\gamma$ is at least $\gamma^\plus$-supercompact.
Please note that this terminology is somewhat strict, with mere
measurability being insufficient. The {\df lottery sum $\oplus
{\cal A}$} of a collection ${\cal A}$ of partial orderings,
defined in \cite{LotteryPreparation}, is the set $\{\langle \Q, q
\rangle : \Q \in {\cal A}$ and $q \in \Q\} \cup \{{\one}\}$,
ordered with $\one$ above everything and $\langle \Q, q \rangle
\le \langle \Q', q' \rangle$ if and only if $\Q = \Q'$ and $q \le
q'$ in $\Q$. Intuitively, a generic object for $\oplus {\cal A}$
selects a ``winning'' poset from ${\cal A}$ and then forces with
it. A forcing notion $\Q$ is {\df stratified} when for any
regular cardinal $\eta$ in the extension the forcing $\Q$ factors
in the ground model as $\Q_0*\Q_1$, in the sense of having
isomorphic complete Boolean algebras, where $\card{\Q_0}\leq\eta$
and $\forces_{\Q_0} \Q_1$ is $\leqeta$-distributive. It follows,
as we observe below, that $\Q_0$ is also stratified. A {\df
non-overlapping iteration} is a forcing iteration $\P$ where the
forcing $\Q_\gamma$ at any stage $\gamma$ in $\P$ is
${\leq}\card{\P_\beta}$-strategically closed for any
$\beta<\gamma$. A forcing notion $\Q$ {\df admits a gap at
$\delta$} if it can be factored as $\Q_0 * \Q_1$, where
$\card{\Q_0} < \delta$ and $\forces_{\Q_0} \Q_1$ is
${\le}\delta$-strategically closed. We denote by $\Cof_\kappa$ the
class of ordinals of cofinality $\kappa$, by
$\add(\theta_1,\theta_2)$ the canonical forcing that adds
$\theta_2$ many Cohen subsets to $\theta_1$ with conditions of
size less than $\theta_1$ and by $\coll(\theta_1,\theta_2)$ the
canonical forcing to collapse $\theta_2$ to $\theta_1$.

\Observation. If\/ $\Q_0*\Q_1$ witnesses the stratification of a
stratified poset $\Q$ at $\eta$, then $\Q_0$ is itself
stratified.\label{Stratified}

\Proof: Since $\card{\Q_0}\leq\eta$, it suffices to stratify
$\Q_0$ only at regular $\zeta<\eta$. Let $\R_0*\R_1$ be the
stratification of $\Q$ at $\zeta$, and suppose
$V[G_0*G_1]=V[H_0*H_1]$ exhibits the equivalent representations
of the forcing extension by $\Q_0*\Q_1$ or $\R_0*\R_1$,
respectively. Since $H_0\in V[G_0][G_1]$ and the $G_1$ forcing is
$\leqeta$-distributive, it follows that $H_0\in V[G_0]$ and so
$V[G_0]=V[H_0][G_0/ H_0]$ for some (quotient) forcing generic
$G_0/ H_0\of\Q_0/ H_0$. That is, $\Q_0$ factors as $\R_0*(\Q_0/
H_0)$. And since $V[H_0]\of V[H_0][G_0/ H_0]\of V[H_0][H_1]$ and
$H_1$ adds no $\zeta$-sequences over $V[H_0]$, it follows that
$G_0/ H_0$ also adds no $\zeta$-sequences over $V[H_0]$. So the
quotient forcing is $\leqzeta$-distributive and we have witnessed
the stratification of $\Q_0$ at $\zeta$.\QED

The careful reader will observe that the definition we gave above
for the level-by-level agreement between strong compactness and
supercompactness presents two distinct departures from a full
general agreement between $\eta$-strong compactness and
$\eta$-supercompactness for every $\eta$. These departures omit
the cases from such a level of agreement that are known to be
generally impossible. The first departure is that we only demand
agreement between $\eta$-strong compactness and
$\eta$-supercompactness when $\eta$ is regular. The reason for
doing so is that Magidor has proved (see \cite[Lemma
7]{ApterShelah495}) that if $\kappa$ is supercompact and $\eta$
is the least strong limit cardinal above $\kappa$ of cofinality
$\kappa$, then there is an $\eta$-supercompactness embedding
$j:V\to M$ such that $\kappa$ is $\eta$-strongly compact but not
$\eta$-supercompact in $M$; consequently, such counterexamples
exist unboundedly often below $\kappa$ as well. The argument also
works for singular $\eta$ of arbitrary cofinality above $\kappa$,
the basic point being that if $\gamma$ is $\lteta$-strongly
compact and $\eta$ is singular with cofinality at least $\gamma$,
then $\gamma$ is $\eta$-strongly compact (but needn't be
$\eta$-supercompact and the least such $\gamma$ cannot be
$\eta$-supercompact). For singular $\eta$ of cofinality less than
$\gamma$ the question is moot because any cardinal $\gamma$ is
$\eta$-strongly compact or $\eta$-supercompact if and only if it
is $\eta^{\ltgamma}$-strongly compact or
$\eta^{\ltgamma}$-supercompact, respectively, rising to the case
of $\eta^{\ltgamma}$, which for such $\eta$ is at least
$\eta^\plus$. And so we restrict our attention to regular degrees
of compactness.

The second restriction is to ignore the case when $\gamma$ is a
measurable limit of cardinals that are $\eta$-strongly compact.
We do this because Menas \cite{Menas74} has shown that such
cardinals are necessarily $\eta$-strongly compact, but not
necessarily $\eta$-supercompact; indeed, if $\eta \ge 2^\gamma$,
then the least such $\gamma$ cannot be $\eta$-supercompact.
Historically, this is how Menas first showed that the notions of
a strongly compact cardinal and a supercompact cardinal are not
identical: the least measurable limit of strongly compact
cardinals will be strongly compact but not supercompact. Later,
this was improved by Magidor \cite{Magidor76} to show that in
fact the least measurable cardinal can be strongly compact.

In summary, the two restrictions in the definition of
level-by-level agreement between strong compactness and
supercompactness omit exactly the cases where such an agreement
is known to be impossible. The point and main contribution of
\cite{ApterShelah495} is that in all other cases, agreement is
possible. In truth, however, the second restriction on
level-by-level agreement will not arise in this article, and can
be safely ignored, because in all the models in which we obtain
level-by-level agreement here, there will be no measurable limits
of cardinals with the same degree of strong compactness.

Finally, we would like to point out that a certain amount of
level-by-level agreement comes for free, namely, a cardinal
$\gamma$ is $\gamma$-strongly compact if and only if it is
$\gamma$-supercompact, since these both are equivalent to
measurability. Therefore, since also $\eta$-supercompactness
directly implies $\eta$-strong compactness for any cardinal, in
order to prove the level-by-level agreement between strong
compactness and supercompactness, it suffices to show that any
cardinal $\gamma$ that is $\eta$-strongly compact for a regular
cardinal $\eta>\gamma$ is also $\eta$-supercompact.

\Section A Surprising Incompatibility \label{Incompatibility}

We begin by proving that if there is a sufficiently large
cardinal above the supercompact cardinal, then the answer to
Question \ref{Open Question} is {\it no}. This result therefore
identifies a surprising incompatibility, a tension between
indestructibility and the level-by-level agreement of strong
compactness and supercompactness in the presence of too many
large cardinals.

\Theorem Incompatibility Theorem. If\/ $\kappa$ is an
indestructible supercompact cardinal and there is a level-by-level
agreement between strong compactness and supercompactness below
$\kappa$, then no cardinal $\lambda$ above $\kappa$ is
$2^\lambda$-supercompact.\label{No}

Theorem \ref{No} is due to the first author and was established
by him in September 1999 during a trip to Japan. To prove this
theorem, we will show that by $\ltkappa$-directed closed forcing
we can force $\lambda$ to violate the level-by-level agreement
between strong compactness and supercompactness. Since the
supercompactness of $\kappa$ will be preserved, it follows by a
simple reflection argument that there must be unboundedly many
violations of the level-by-level agreement below $\kappa$; this
contradicts the fact that the level-by-level agreement below
$\kappa$ is preserved to the forcing extension.

So let us begin with the following:

\SubLemma. If a cardinal $\lambda$ is $2^\lambda$-supercompact,
then there is a forcing extension $V^\P$ in which $\lambda$ is
$\lambda^\plus$-strongly compact but not
$\lambda^\plus$-supercompact. In $V^\P$, one can arrange that
$2^\lambda=\lambda^\plus$ and $\lambda$ has, as a measurable
cardinal, trivial Mitchell rank. Consequently, $\lambda$ will not
be even $(\lambda+2)$-strong in $V^\P$. What's more, for any
$\delta<\lambda$, the forcing $\P$ can be chosen to be
$\leqdelta$-directed closed and with a gap below
$\lambda$.\label{StrCnotSC}

\Proof: Standard arguments establish that if $\lambda$ is
$2^\lambda$-supercompact, then this can be preserved to a forcing
extension in which $2^\lambda=\lambda^\plus$ (one simply forces
$2^\gamma=\gamma^\plus$ with $\add(\gamma^\plus,1)$ at
sufficiently many stages $\gamma\leq\lambda$ in a reverse Easton
iteration). This iteration admits a gap between any two
nontrivial stages of forcing, and by starting the iteration
beyond any particular $\delta<\lambda$, we may ensure that it is
$\leqdelta$-closed. Afterwards, we may directly force
$2^{\lambda^\plus}=\lambda^\plusplus$ by adding a Cohen subset to
$\lambda^\plusplus$; since this adds no subsets to
$P_\lambda\lambda^\plus$, it therefore preserves the
$\lambda^\plus$-supercompactness of $\lambda$. So let us assume
without loss of generality that we have already performed this
forcing, if necessary, and that $2^\lambda=\lambda^\plus$ and
$2^{\lambda^\plus}=\lambda^\plusplus$ in $V$.

Let $\P$ be the reverse Easton $\lambda$-iteration which at stage
$\gamma$ forces with $\Q_\gamma=\add(\gamma,1)$, provided that
$\gamma$ is above $\delta$ and measurable in $V$. This forcing is
$\leqdelta$-closed and admits a gap between any two nontrivial
stages of forcing. Suppose that $G\of\P$ is $V$-generic.

We claim that $\lambda$ has trivial Mitchell rank in $V[G]$. If
not, then there would be an embedding $j:V[G]\to M[j(G)]$ with
critical point $\lambda$ for which $M[j(G)]$ is closed under
$\lambda$-sequences in $V[G]$ and $\lambda$ is measurable in
$M[j(G)]$. By the Gap Forcing Theorem of \cite{GapForcing}
applied in $M[j(G)]$, it follows that $\lambda$ is measurable in
$M$ and consequently a stage of nontrivial forcing in $j(\P)$.
Factoring $j(\P)$ as $\P*\add(\lambda,1)*\Ptail$, it follows that
$j(G)$ must be $G*A*\Gtail$ for some $M[G]$-generic Cohen subset
$A\of\lambda$. Since every subset of $\lambda$ in $V$ is in $M$
and the forcing $\P$ has size $\lambda$, it follows that every
subset of $\lambda$ in $V[G]$ is in $M[G]$. In particular, every
dense subset of $\add(\lambda,1)^{V[G]}=\add(\lambda,1)^{M[G]}$
from $V[G]$ is in $M[G]$ and so $A$ is actually $V[G]$-generic as
well. Since this contradicts the fact that $A\in V[G]$, there can
be no such embedding $j$ and so $\lambda$ has trivial Mitchell
rank in $V[G]$. It follows that $\lambda$ is neither
$(\lambda+2)$-strong nor $\lambda^\plus$-supercompact in $V[G]$.

We claim nevertheless that $\lambda$ remains
$\lambda^\plus$-strongly compact in $V[G]$. For this, we use a
technique of Magidor, unpublished by him but exposited in
\cite{ApterCummings2000:IdentityCrises}, \cite{AptCum2},
\cite{ApterA}, \cite{ApterB}, \cite{ApterC}, and \cite{ApterD}.
Let $j_0:V\to M$ be a $\lambda^\plus$-supercompactness embedding
generated by a normal fine measure on $P_\lambda\lambda^\plus$
and $h:M\to N$ an ultrapower embedding by a measure on $\lambda$
of minimal Mitchell rank in $M$, so that $\lambda$ is not
measurable in $N$. Let $j=h\compose j_0$ be the combined
embedding; it witnesses the $\lambda^\plus$-strong compactness of
$\lambda$. We will lift this embedding to $j:V[G]\to N[j(G)]$ so
as to witness the $\lambda^\plus$-strong compactness of $\lambda$
in $V[G]$.

Consider the forcing $j(\P)$, factored as
$\P*\P_{\lambda,h(\lambda)}*\add(h(\lambda),1)*\P_{h(\lambda),j(\lambda)}$,
and the forcing $j_0(\P)$, factored as
$\P*\add(\lambda,1)*\Ptail$. With this notation, for example,
$h(\P*\add(\lambda,1))=\P*\P_{\lambda,h(\lambda)}*\add(h(\lambda),1)$.

Now let $\Pterm$ be the term forcing poset for $\Ptail$ over
$\P*\add(\lambda,1)$ (see \cite{Foreman83} for the first published
account of term forcing, or \cite[1.2.5, p. 8]{Cummings92GCH}; the
notion is originally due to Richard Laver). That is, $\Pterm$
consists of (sufficiently many) $\P*\add(\lambda,1)$-names for
elements of $\Ptail$, ordered by $\tau\leq\sigma$ if and only if
$\one\forces\tau\leq\sigma$. As in the proof of \cite[Lemma
3.2]{ApterD}, a full collection of names, meaning that any name
forced by $\one$ to be in $\Ptail$ is forced by $\one$ to be
equal to one of them, can be found of size $j(\lambda)$ in $M$,
which has cardinality $\lambda^\plus$ in $V$. Further, since
$\Ptail$ is forced to be $\leqlambda^\plus$-directed closed, it
is easy to see that $\Pterm$ is $\leqlambda^\plus$-directed
closed in $M$, and hence also $\leqlambda^\plus$-directed closed
in $V$. Since $M$ has only $j(\lambda^\plus)$ many dense sets for
$\Pterm$, and this has cardinality $\lambda^\plusplus$ in $V$, we
may by the usual diagonalization techniques (see, e.g.
\cite{FragileMeasurability}) construct an $M$-generic filter
$\Gterm\of\Pterm$ in $V$. And since $h$ is the ultrapower by a
measure on $\lambda$ and $\Pterm$ is $\leqlambda$-closed, it
follows that $h\image\Gterm$ is $N$-generic for the term forcing
for $\P_{h(\lambda)+1,j(\lambda)}$ over $\P_{h(\lambda)+1}$ (see
\cite{Foreman83} or  \cite[1.2.2, Fact 2]{Cummings92GCH}). Thus,
we may lift the embedding $h$ to $h:M[\Gterm]\to N[h(\Gterm)]$.
And since $G\of\P$ is $V$-generic, it is also
$N[h(\Gterm)]$-generic, and so we may form the extension
$N[h(\Gterm)][G]$.

Since $\lambda$ is not measurable in $N$, the stage $\lambda$
forcing in $j(\P)$ is trivial, and so $\P_{\lambda,h(\lambda)}$ is
$\leqlambda$-closed in $N[G]$. Since the $h(\Pterm)$ forcing is
highly closed in $N$, it cannot affect closure down at $\lambda$,
so the forcing $\P_{\lambda,h(\lambda)}$ is $\leqlambda$-closed
in $N[h(\Gterm)][G]$. Going one step more, the poset
$\P_{\lambda,h(\lambda)}*\add(h(\lambda),1)$ is
$\leqlambda$-closed in $N[h(\Gterm)][G]$, has size $h(\lambda)$
there, which has size $\lambda^\plus$ in $V$, and
$N[h(\Gterm)][G]$ is closed under $\lambda$-sequences in $V[G]$.
Thus, by the usual diagonalization techniques, we may construct in
$V[G]$ (actually, we can do it in $M[G]$) an
$N[h(\Gterm)][G]$-generic $G_{\lambda,h(\lambda)+1}$ for this
much of $j(\P)$.

By the fundamental property of term forcing (see \cite{Foreman83}
or \cite[1.2.5, Fact 1]{Cummings92GCH}), in
$N[h(\Gterm)][G_{h(\lambda)+1}]$ we may construct an
$N[G_{h(\lambda)+1}]$-generic filter
$G_{h(\lambda)+1,j(\lambda)}$ from $h(\Gterm)$. Putting these
filters together, let
$j(G)=G_{h(\lambda)+1}*G_{h(\lambda)+1,j(\lambda)}$ and lift the
embedding to $j:V[G]\to N[j(G)]$ in $V[G]$. The attentive reader
will observe that we used the term forcing only to help construct
$G_{h(\lambda)+1,j(\lambda)}$. Now that we have done so, we may
discard $\Gterm$ and $h(\Gterm)$.

To see that this lifted embedding witnesses the
$\lambda^\plus$-strong compactness of $\lambda$ in $V[G]$, let
$s=h(j_0\image\lambda^\plus)$. One can now easily check that $s\in
N$, $s\of j(\lambda^\plus)$,
$\card{s}=h(\lambda^\plus)<j(\lambda)$ and
$j\image\lambda^\plus\of s$. Thus, $s$ induces a fine measure on
$P_\lambda\lambda^\plus$ in $V[G]$, as desired.\footnote{In fact,
one can show more: the lifted embedding $j$ itself is the
ultrapower by a fine measure on $P_\lambda\lambda^\plus$, rather
than merely inducing such a measure; one need only modify $s$ by
adding an ordinal on top from which we can definably recover the
ordinal seed for the measure used in the ultrapower embedding
$h$. The resulting seed $s^*$ will still induce a fine measure
$\mu$ on $P_\lambda\lambda^\plus$, but since the hull of $s^*$
with $\ran(j)$ is all of $M[j(G)]$, the ultrapower by $\mu$ will
be precisely $j$. We refer readers to \cite{Seeds} for an account
of this seed technique and terminology.}\QED

It remains to check one last detail before the proof of Theorem
\ref{No} is complete: we need to know that $\lambda$ is not one
of the exceptional cases excluded from our definition of
level-by-level agreement. Let us assume without loss of
generality that we worked with the least $\lambda$ above $\kappa$
that was $2^\lambda$-supercompact. In this case, because of the
level-by-level agreement, $\lambda$ cannot be in $V$ a measurable
limit of cardinals that are $2^\lambda$-strongly compact. By the
Gap Forcing Theorem \cite{GapForcing}, the forcing iteration does
not increase the degree of strong compactness of any cardinal, so
in the extension $\lambda$ is $\lambda^\plus$-strongly compact,
but not $\lambda^\plus$-supercompact and not a measurable limit
of cardinals that are $\lambda^\plus$-strongly compact. Thus, it
is truly a violation of the level-by-level agreement between
strong compactness and supercompactness. By reflecting this
violation below $\kappa$, we contradict our assumption that there
was a level-by-level agreement there, and the proof of Theorem
\ref{No} is complete.

The initial proof of Lemma \ref{StrCnotSC} used the original form
of Magidor's method, which iteratively adds stationary
non-reflecting sets; later, we saw that an appeal to the Gap
Forcing Theorem allowed us to simplify this to just iterated
Cohen forcing.

Before concluding this section, we would like to call attention
to the fact that our proof of Theorem \ref{No} does not fully use
the hypothesis that $\kappa$ is an indestructible supercompact
cardinal. First, one can easily check that the proof uses only
that the $2^{2^\lambda}$-supercompactness of $\kappa$ is
indestructible. But in fact, the reflection argument that we used
to bring the violation of the level-by-level agreement at
$\lambda$ down below $\kappa$ does not actually rely on the
supercompactness of $\kappa$ at all, but only on its strongness.
Therefore, we have actually proved the following theorem.

\Theorem. If\/ $\kappa$ is an indestructible strong cardinal and
there is a level-by-level agreement between strong compactness and
supercompactness below $\kappa$, then no cardinal $\lambda$ above
$\kappa$ is $2^\lambda$-supercompact. Indeed, $\kappa$ need only
be indestructibly $(\lambda+2)$-strong.

What's more, since the iteration up to $\lambda$ can be arranged
to be $\leqkappa$-closed, we only need the strongness of $\kappa$
to be indestructible by $\leqkappa$-closed forcing (or even less:
indestructible by $\leqkappa^\plus$-closed forcing, etc.). This
is interesting because \cite{GitikShelah89} shows that any strong
cardinal $\kappa$ can be made indestructible by
$\leqkappa$-closed forcing. So our argument shows that this
amount of indestructibility for a strong cardinal is incompatible
with a level-by-level agreement between strong compactness and
supercompactness if there are large enough cardinals above.

But actually, we don't even need $\kappa$ to be strong for the
reflection to work. Since the violation of the level-by-level
agreement at $\lambda$ is witnessed in $V_{\lambda+2}$, a rank
initial segment of $V$, it is enough if $V_\kappa\elesub_2 V$,
that is, if $\kappa$ is {\df $\Sigma_2$-reflecting}. Therefore, we
have proved:

\Theorem. If\/ $\kappa$ is an indestructible $\Sigma_2$-reflecting
cardinal and there is a level-by-level agreement between strong
compactness and supercompactness below $\kappa$, then no cardinal
$\lambda$ above $\kappa$ is $2^\lambda$-supercompact.

Indeed, it is enough if the relation $V_\kappa\elesub_2 V_\theta$
is indestructible by $\ltkappa$-directed closed forcing in
$V_\theta$ for some $\theta>\lambda+1$ (or even indestructible
just by $\leqkappa$-directed closed forcing, etc.).

\Section Level-by-level agreement with near indestructibility
\label{NearSec1}

In this section we provide an affirmative answer to a weakened
form of Question \ref{Open Question}, showing that a supercompact
cardinal can be nearly indestructible in the presence of a
level-by-level agreement between strong compactness and
supercompactness. Specifically, we will obtain the level-by-level
agreement in a model with a supercompact cardinal $\kappa$ that
is indestructible by all stratified $\ltkappa$-directed closed
forcing and more. Recall our definition that a poset is {\df
stratified} when for every regular cardinal $\eta$ in the
extension it factors in the ground model as $\Q_0*\Q_1$, in the
sense of having isomorphic complete Boolean algebras, where
$\card{\Q_0}\leq\eta$ and $\forces_{\Q_0}\Q_1$ is
$\leqeta$-distributive. It follows, as we proved in Observation
\ref{Stratified}, that $\Q_0$ is also stratified. We say that
$\Q$ is {\df stratified above $\kappa$} if such a factorization
exists for $\eta$ above $\kappa$. Numerous examples of stratified
forcing, including Cohen forcing and collapsing posets as well as
iterations of these and many others, are indicated in Corollary
\ref{Strat}.

The proof will proceed by an iteration that we call the lottery
preparation preserving level-by-level agreement, and we regard
the resulting model as currently the most natural candidate for an
affirmative answer to Question \ref{Open Question}, if any exists.
In particular, while we have been able to prove so far only that
the supercompactness of $\kappa$ is indestructible there by any
stratified $\ltkappa$-directed closed forcing and a few others, we
know of no obstacle to it being fully indestructible there by all
$\ltkappa$-directed closed forcing.

\Theorem. Suppose that $\kappa$ is supercompact and no cardinal
is supercompact up to a cardinal $\lambda$ which is itself\/
$2^\lambda$-supercompact. Then there is a forcing extension
satisfying a level-by-level agreement between strong compactness
and supercompactness in which $\kappa$ remains supercompact and
becomes indestructible by any stratified $\ltkappa$-directed
closed forcing and more. Indeed, the supercompactness of $\kappa$
becomes indestructible by any $\ltkappa$-directed closed forcing
that is stratified above $\kappa$.\label{Near}

\Proof: We may assume without loss of generality, by forcing if
necessary, that the \GCH\ holds and further, by forcing with the
notion in \cite{ApterShelah495} if necessary, that in $V$ there
is already a level-by-level agreement between strong compactness
and supercompactness. Since these forcing notions admit a very
low gap, by the Gap Forcing Theorem \cite{GapForcing} they do not
increase the degree of supercompactness or (since the forcing
notions are also mild in the sense of \cite{GapForcing}) strong
compactness of any cardinal. It follows that no cardinal in $V$ is
supercompact up to a partially supercompact cardinal. Let $\P$ be
the reverse Easton support $\kappa$-iteration which begins by
adding a Cohen real and then has nontrivial forcing only at later
stages $\gamma$ that are inaccessible limits of partially
supercompact cardinals. At such a stage $\gamma$ in $\P$,
assuming $V^{\P_\gamma}$ has a level-by-level agreement between
strong compactness and supercompactness, we force with the
lottery sum of all $\ltgamma$-directed closed posets $\Q$, of
size less than the next partially supercompact cardinal, that
preserve this level-by-level agreement. (Please note that this
will always include the trivial poset.) The iteration $\P$ is an
example of a modified lottery preparation of the type used in
\cite{AsYouLikeIt}, \cite{LotteryPreparation}, and \cite{ApterC}.
Supposing $G\of\P$ is $V$-generic, we will refer to the iteration
$\P$ and the resulting model $V[G]$ as the lottery preparation
preserving level-by-level agreement.\footnote{And whenever we use
this terminology, we implicitly assume that the ground model $V$
has a level-by-level agreement between strong compactness and
supercompactness, as well as the \GCH, and that no cardinal in
$V$ is supercompact up to a partially supercompact cardinal.}

\SubLemma. If\/ $\gamma$ is $\eta$-supercompact in $V$ for a
regular cardinal $\eta > \gamma$, then this remains true in
$V[G_\gamma]$. In particular, after the lottery preparation for
preserving level-by-level agreement $V[G]$, the cardinal $\kappa$
remains fully supercompact.\label{SubA}

\Proof: (This same observation was essentially made in
\cite[Lemma 2.1]{AsYouLikeIt}, the modified lottery preparations
resembling as they do the partial Laver preparations.) We may
assume that $\gamma$ is a limit of partially supercompact
cardinals, since otherwise the forcing $\P_\gamma$ is equivalent
to forcing that is small with respect to $\gamma$ and the result
is immediate by \cite{LevySolovay67}. Let $j:V\to M$ be an
$\eta$-supercompactness embedding with critical point $\gamma$
such that $\gamma$ is not $\eta$-supercompact in $M$. Since
$\gamma$ is $\lteta$-supercompact in $M$ and no cardinal is
supercompact up to a partially supercompact cardinal, it follows
that the next partially supercompact cardinal above $\gamma$ in
$M$ is at least $\eta$; further, $\eta$ itself is not even
measurable in $M$ because if it were, then the
$\lteta$-supercompactness of $\gamma$ would imply that $\gamma$
is $\eta$-supercompact in $M$, contrary to our assumption. In
short, the next partially supercompact cardinal in $M$ above
$\gamma$ is strictly above $\eta$. Thus, below a condition opting
for trivial forcing at stage $\gamma$ in $j(\P_\gamma)$, we may
factor $j(\P_\gamma)$ as $\P_\gamma*\Ptail$, where $\Ptail$ is
$\leqeta$-closed. Thus, by the usual diagonalization techniques,
we may construct in $V[G_\gamma]$ an $M[G_\gamma]$-generic filter
$\Gtail\of\Ptail$ and lift the embedding to $j:V[G_\gamma]\to
M[j(G_\gamma)]$ with $j(G_\gamma)=G_\gamma*\Gtail$. This
embedding witnesses that $\gamma$ is $\eta$-supercompact in
$V[G_\gamma]$, as desired.\QED

\SubLemma. For any ordinal\/ $\gamma$, there is a level-by-level
agreement between strong compactness and supercompactness in
$V[G_\gamma]$. In particular, the lottery preparation preserving
level-by-level agreement $V[G]$ really does preserve the
level-by-level agreement between strong compactness and
supercompactness.

\Proof: Suppose inductively that the result holds below $\gamma$
and consider $V[G_\gamma]$. Since successor stages of forcing
always preserve the level-by-level agreement if it exists, we may
assume that $\gamma$ is a limit of stages of forcing, and hence a
limit of partially supercompact cardinals. For any
$\delta<\gamma$, our induction hypothesis guarantees a
level-by-level agreement for $\delta$ in $V[G_{\delta+1}]$, and
by the Gap Forcing Theorem [Ham], $\delta$ is neither strongly
compact nor supercompact up to a partially supercompact cardinal
in that model. Therefore, since the later non-trivial stages of
forcing are closed beyond the next inaccessible limit of
partially supercompact cardinals, which by our assumptions and
the level-by-level agreement is beyond the degree of strong
compactness or supercompactness of $\delta$, it follows that
there is a level-by-level agreement for $\delta$ in
$V[G_\gamma]$. By \cite{LevySolovay67}, cardinals above $\gamma$
are not affected by small forcing or anything equivalent to small
forcing, and so we need only consider the cardinal $\gamma$
itself. Suppose accordingly that $\gamma$ is $\eta$-strongly
compact in $V[G_\gamma]$ for a regular cardinal $\eta>\gamma$; we
will show it is also $\eta$-supercompact there. By the Gap
Forcing Theorem \cite{GapForcing}, we know that $\gamma$ is
$\eta$-strongly compact in $V$, and hence by the level-by-level
agreement, it is also $\eta$-supercompact there. So by the
previous lemma $\gamma$ is $\eta$-supercompact in $V[G_\gamma]$,
as desired.\QED

\SubLemma. If\/ $\gamma$ is an inaccessible limit of partially
supercompact cardinals, then any stratified $\ltgamma$-directed
closed forcing $\Q$ over $V[G_\gamma]$ preserves the
level-by-level agreement between strong compactness and
supercompactness for $\gamma$ and smaller cardinals. Indeed, $\Q$
need only be stratified at regular cardinals $\eta$ above
$\gamma$.\label{SubB}

\Proof: Suppose that the result holds for cardinals below $\gamma$
(with full Boolean value) and that $\Q$ is $\ltgamma$-directed
closed in $V[G_\gamma]$ and stratified for regular $\eta$ above
$\gamma$. From the closure of $\Q$ and the fact that no cardinal
is supercompact beyond a partially supercompact cardinal in $V$
and hence also (by the Gap Forcing Theorem \cite{GapForcing}) in
$V[G_\gamma]$, one sees that it must preserve the level-by-level
agreement between strong compactness and supercompactness for all
cardinals below $\gamma$. So we consider the cardinal $\gamma$
itself.

Suppose that $\gamma$ is $\eta$-strongly compact in
$V[G_\gamma]^\Q$ for some regular cardinal $\eta>\gamma$. We will
show that $\gamma$ is $\eta$-supercompact there as well. By the
Gap Forcing Theorem \cite{GapForcing}, we know that $\gamma$ is
$\eta$-strongly compact in $V$, and hence by the level-by-level
agreement, it is $\eta$-supercompact there as well. Fix an
$\eta$-supercompactness embedding $j:V\to M$ with $\gamma$ not
$\eta$-supercompact in $M$. As in Lemma \ref{SubA} it follows that
the next partially supercompact cardinal of $M$ above $\gamma$ is
above $\eta$. Since $\Q$ is stratified above $\gamma$, we may
factor $\Q$ as $\Q_0*\Q_1$, where $\card{\Q_0}\leq\eta$ and
$\forces_{\Q_0}\Q_1$ is $\leqeta$-distributive. By cardinality
considerations, forcing with $\Q_1$ adds no new subsets of
$P_\gamma\eta$ over $V[G_\gamma]^{\Q_0}$, and so it suffices for
us to show that $\gamma$ is $\eta$-supercompact in
$V[G_\gamma][g]$, where $g\of \Q_0$ is $V[G_\gamma]$-generic. The
cardinal $\gamma$ is an inaccessible limit of partially
supercompact cardinals in $M$, and so there is a lottery at stage
$\gamma$ in $j(\P_\gamma)$. Furthermore, since the induction
hypothesis holds up to $j(\gamma)$ in $M$ and $\Q_0$ is
$\ltgamma$-directed closed, has size at most $\eta$ in $M$ and,
by Observation \ref{Stratified}, is stratified above $\gamma$
(and this can be seen in $M[G_\gamma]$), it follows that $\Q_0$
is allowed to appear in the stage $\gamma$ lottery of
$j(\P_\gamma)$. Thus, below a condition opting for $\Q_0$ in this
lottery, we may factor the forcing $j(\P_\gamma)$ as
$\P_\gamma*\Q_0*\Ptail$, where $\Ptail$ is $\leqeta$-closed. Now
the usual diagonalization arguments apply and we may construct in
$V[G_\gamma][g]$ an $M[G_\gamma][g]$-generic filter
$\Gtail\of\Ptail$ and lift the embedding to $j:V[G_\gamma]\to
M[j(G_\gamma)]$ where $j(G_\gamma)=G_\gamma*g*\Gtail$. Since
$j\image g\of j(\Q_0)$ is a directed subset of cardinality less
than $j(\gamma)$ in $M[j(G_\gamma)]$, we may find a (master)
condition $p$ below it. Working now below $p$ we diagonalize to
construct an $M[j(G_\gamma)]$-generic object $j(g)\of j(\Q_0)$
and lift the embedding fully to $j:V[G_\gamma][g]\to
M[j(G_\gamma)][j(g)]$, thereby witnessing the
$\eta$-supercompactness of $\gamma$ in $V[G_\gamma][g]$, as
desired.\QED

\SubLemma. After the lottery preparation preserving
level-by-level agreement $V[G]$, the supercompactness of\/
$\kappa$ becomes indestructible by any $\ltkappa$-directed closed
forcing that is stratified above $\kappa$.

\Proof: Suppose that $\Q$ is stratified above $\kappa$ and
$\ltkappa$-directed closed in $V[G]$, and assume $g\of\Q$ is
$V[G]$-generic. Select any regular $\eta\geq\card{\Q}$ and an
$\eta$-supercompactness embedding $j:V\to M$ for which $\kappa$
is not $\eta$-supercompact in $M$. It follows as in Lemma
\ref{SubB} that $\Q$ is allowed in the stage $\kappa$ lottery and
the diagonalization arguments given in Lemma \ref{SubB} allow us
to lift the embedding to $j:V[G][g]\to M[j(G)][j(g)]$. This
witnesses the $\eta$-supercompactness of $\kappa$ in $V[G][g]$, as
desired.\QED

This completes the proof of the theorem.\QED

\Corollary. After the lottery preparation preserving
level-by-level agreement, the supercompactness of\/ $\kappa$
becomes indestructible by:
\begin{enumerate}
\item the Cohen forcing $\add(\kappa,1)$,
\item the Cohen forcing $\add(\theta,1)$ for any
regular $\theta$ above $\kappa$,
\item the Cohen forcing $\add(\kappa,\kappa^\plus)$,
\item indeed, any $\ltkappa$-directed closed forcing of size at most
$\kappa^\plus$,
\item the collapse forcing $\coll(\kappa,\theta)$ for any $\theta$ above
$\kappa$,
\item the collapse forcing $\coll(\theta_1,\theta_2)$ whenever
$\kappa\leq\theta_1\leq\theta_2$ and $\theta_1$ is regular,
\item the L\'evy collapse forcing $\coll(\theta_1, {<}\theta_2)$
whenever $\kappa \le \theta_1 \leq \theta_2$ and $\theta_1$ is
regular,
\item the forcing to add a stationary non-reflecting subset
of\/ $\theta\intersect\Cof_\kappa$ for any regular $\theta>
\kappa$,
\item and any non-overlapping reverse Easton support iteration of these forcing
notions.
\end{enumerate}
\label{Strat}

\Proof: These are all $\ltkappa$-directed closed and either
completely stratified or stratified above $\kappa$. Indeed, since
$\add(\theta,1)=\coll(\theta,\theta)$, we can see that cases 1, 2
and 5 are special cases of $\coll(\theta_1,\theta_2)$, in case 6.
And it is easy to see that this and $\coll(\theta_1,
{<}\theta_2)$, for $\kappa \leq \theta_1\leq\theta_2$ and
$\theta_1$ regular, are stratified: any regular cardinal $\eta$ in
the extension must be either at most $\theta_1$ or at least
$\theta_2$, and one can trivially factor the forcing. (Note that
we use the \GCH\ in $V$ and the fact $\card{\P} = \kappa$ in
order to know that
$\card{\coll(\theta_1,\theta_2)}=\theta_2^{\lttheta_1}\leq\eta$ if
$\eta$ is regular and at least $\theta_2$.) The cases of
$\add(\kappa,\kappa^\plus)$ or any $\ltkappa$-directed closed
$\Q$ of size at most $\kappa^\plus$ are clearly stratified above
$\kappa$. The conditions of the forcing to add a stationary
non-reflecting subset to $\theta\intersect\Cof_\kappa$ are simply
the bounded subsets $s\of\theta\intersect\Cof_\kappa$ which are
not stationary in their supremum nor have any initial segment
stationary in its supremum, ordered by end-extension. This
forcing is $\ltkappa$-directed closed by simply taking unions of
conditions. To see that it is stratified, we note that for any
$\eta<\theta$, the forcing is $\leqeta$-strategically closed and
hence $\leqeta$-distributive, and for any regular
$\eta\geq\theta$, by \GCH\ in $V$ and the fact $\card{\P} =
\kappa$, the forcing has size less than or equal to $\eta$. So in
each case the factorization is trivial. Finally, it is easy to
see that a non-overlapping iteration of these posets remains
stratified above $\kappa$.\QED

We would like to call attention to a sense in which the
hypothesis in Theorem \ref{Near}, namely, that no cardinal is
supercompact up to a larger cardinal $\lambda$ that is itself
$2^\lambda$-supercompact, is optimal in light of Theorem \ref{No}
and Corollary \ref{Strat}. Specifically, Theorem \ref{No}
establishes that if one has a level-by-level agreement and a
supercompact cardinal $\kappa$ that is indestructible by iterated
Cohen forcing, then no larger cardinal $\lambda$ is
$2^\lambda$-supercompact. Conversely, if there is no such
$\lambda$ above $\kappa$, then Corollary \ref{Strat} shows that
this amount of indestructibility is possible. In this sense, the
conclusion of Theorem \ref{No} tightly matches the hypothesis of
Theorem \ref{Near}.

The next theorem, however, shows that in another sense Theorem
\ref{Near} and Corollary \ref{Strat} are not optimal; there is
more indestructibility than we have claimed so far.  Specifically,
while each of the posets mentioned in Corollary \ref{Strat}
preserves the \GCH, and indeed, a straightforward factor argument
shows that any stratified forcing whatsoever preserves the \GCH,
our preparatory forcing in fact makes the supercompactness of
$\kappa$ indestructible by $\add(\kappa,\kappa^\plusplus)$, which
of course does not preserve the \GCH.

\Corollary. After the lottery preparation preserving
level-by-level agreement, the supercompactness of $\kappa$
becomes indestructible by\/ $\add(\kappa,
\kappa^{\plusplus})$.\label{gchind}

\Proof: We reiterate our implicit assumption for this preparation
that the ground model $V$ satisfies the \GCH\ and a level-by-level
agreement between strong compactness and supercompactness, and
that no cardinal in $V$ is supercompact up to a partially
supercompact cardinal there.

We begin by showing that for any cardinal $\gamma\leq\kappa$ the
forcing $\add(\gamma,\gamma^\plusplus)$ over $V[G_\gamma]$
preserves the level-by-level agreement between strong compactness
and supercompactness. Suppose inductively that this holds (with
full Boolean value) for all cardinals below $\gamma$. Since
$\add(\gamma,\gamma^\plusplus)$ adds no small sets, it does not
affect the level-by-level agreement between strong compactness and
supercompactness for cardinals below $\gamma$, and since the
forcing is itself small with respect to larger cardinals, it
preserves such agreement above $\gamma$ by \cite{LevySolovay67}.
So it remains only to check the agreement right at $\gamma$.
Accordingly, suppose that $\gamma$ is $\lambda$-strongly compact
for some regular cardinal $\lambda>\gamma$ in $V[G_\gamma][g]$,
where $g\of\add(\gamma,\gamma^\plusplus)$ is
$V[G_\gamma]$-generic; we aim to show that $\gamma$ is also
$\lambda$-supercompact there.

The Gap Forcing Theorem \cite{GapForcing} implies that $\gamma$
is $\lambda$-strongly compact in $V$ and hence also
$\lambda$-supercompact, by the level-by-level agreement there,
witnessed by an embedding $j:V\to M$. We may assume that $\gamma$
is not $\lambda$-supercompact in $M$. We will lift this embedding
to $V[G_\gamma][g]$, thereby witnessing the
$\lambda$-supercompactness of $\gamma$ there.

We note first that $\gamma$ must be a limit of partially
supercompact cardinals, since otherwise the forcing
$\P_\gamma*\add(\gamma,\gamma^\plusplus)$ is equivalent to small
forcing followed by $\ltgamma$-closed forcing that adds a subset
to $\gamma$; but by \cite{SmallForcing}, all such forcing destroys
the measurability of $\gamma$, contrary to our assumption that
$\gamma$ is $\lambda$-strongly compact in $V[G_\gamma][g]$. And
since these smaller cardinals are fixed by $j$, it follows now
that $\gamma$ is an inaccessible limit of partially supercompact
cardinals in $M$, and hence a nontrivial stage of forcing in
$j(\P_\gamma)$.

By elementarity the induction hypothesis holds up to $j(\gamma)$
in $M$, and so the forcing $\add(\gamma,\gamma^\plusplus)$, which
is the same in $V[G_\gamma]$ and $M[G_\gamma]$ since
$\gamma^\plus\leq\lambda$, preserves the level-by-level agreement
between strong compactness and supercompactness over
$M[G_\gamma]$. It is therefore allowed to appear in the stage
$\gamma$ lottery of $j(\P_\gamma)$. Below a condition opting for
this poset in that lottery, therefore, the forcing $j(\P_\gamma)$
factors as $\P_\gamma*\add(\gamma,\gamma^\plusplus)*\Ptail$,
where $\Ptail$ is the part of the forcing at stages beyond
$\gamma$. Since the next partially supercompact cardinal in $M$
must be beyond $\lambda$, it follows that $\Ptail$ is
$\leqlambda$-closed in $M[G_\gamma][g]$. Further,
$M[G_\gamma][g]$ is closed under $\lambda$-sequences in
$V[G_\gamma][g]$ and the number of dense subsets of $\Ptail$ in
$M[G_\gamma][g]$ is at most
$\card{j(\gamma^\plus)}^V\leq(\gamma^\plus)^\lambda=\lambda^\plus$.
Therefore, we may simply line up these dense sets in
$V[G_\gamma][g]$, diagonalize to construct an
$M[G_\gamma][g]$-generic filter $\Gtail\of\Ptail$ and lift the
embedding to $j:V[G_\gamma]\to M[j(G_\gamma)]$ with
$j(G_\gamma)=G_\gamma*g*\Gtail$. It remains to lift the embedding
through the forcing $\add(\gamma,\gamma^\plusplus)$.

If $\gamma^\plusplus\leq\lambda$, and this is the easy case, the
usual master condition argument allows us to lift the embedding.
Specifically, since $g\in M[j(G_\gamma)]$ one uses
$j\image\gamma^\plusplus$ to see that $j\image g$ is also in
$M[j(G_\gamma)]$, and since this is a directed subset of
$j(\add(\gamma,\gamma^\plusplus))$ of size less than $j(\gamma)$,
there is a condition $p$, called the master condition, below it.
Diagonalizing below this condition, one builds an
$M[j(G_\gamma)]$-generic filter $g^*\of
j(\add(\gamma,\gamma^\plusplus))$ and lifts the embedding fully
to $j:V[G_\gamma][g]\to M[j(G_\gamma)][j(g)]$, with $j(g)=g^*$,
as desired.

For the only remaining case, the hard case, assume
$\lambda=\gamma^\plus$. In this case $j\image g$ will not be in
$M[j(G_\gamma)]$, and there will be no master condition.
Nevertheless, with care we will still be able to construct a
generic filter extending $j\image g$. (This technique appears in
\cite[p. 119-120]{ApterShelah495} and \cite[p. 555-556]
{Apter99ViolateGCH}, and is similar to a technique involving
strong cardinals in \cite[p. 277--278]{FragileMeasurability}.) We
will construct an $M[j(G_\gamma)]$-generic filter $g^*\of
j(\add(\gamma,\gamma^\plusplus))$ in $V[G_\gamma][g]$ such that
$j\image g\of g^*$, and then lift the embedding fully to
$j:V[G_\gamma][g]\to M[j(G_\gamma)][j(g)]$, with $j(g)=g^*$,
thereby witnessing the $\lambda$-supercompactness of $\gamma$ in
$V[G_\gamma][g]$.

As a first step towards this, suppose that $A\of
j(\add(\gamma,\gamma^\plusplus))$ is a maximal antichain in
$M[j(G_\gamma)]$ and $r\in j(\add(\gamma,\gamma^\plusplus))$ is a
condition that is compatible with every element of $j\image g$. We
will find a condition $r^+\leq r$ that decides $A$ while
remaining compatible with $j\image g$. Since
$\add(\gamma,\gamma^\plusplus)$ is $\gamma^\plus$-c.c., it
follows that the antichain $A$ has size at most $j(\gamma)$ in
$M[j(G_\gamma)]$. Since also the usual supercompactness arguments
show that $j\image\gamma^\plusplus$ is unbounded in
$j(\gamma^\plusplus)$, it follows that $A\of
j(\add(\gamma,\alpha))$ for sufficiently large
$\alpha<\gamma^\plusplus$. Fix such an $\alpha$ such that also
$r\in j(\add(\gamma,\alpha))$ and let $q=\union \bigl(j\image
(g\intersect\add(\gamma,\alpha))\bigr)$. Since $q$ has size at
most $\gamma^\plus$, it is in $M[j(G_\gamma)]$ and is a (master)
condition for $j(\add(\gamma,\alpha))$, which is a complete
subposet of $j(\add(\gamma,\gamma^\plusplus))$. Since $r$ is
compatible with every element of $j\image g$, we know that $r$
and $q$ are compatible. Choose $r^+\in j(\add(\gamma,\alpha))$
below $r$ and $q$ and deciding $A$. We claim that $r^+$ remains
compatible with $j\image g$. To see this, consider any condition
$j(p)$ for $p\in g$. Since $p$ is a partial function from
$\gamma\cross \gamma^\plusplus$ into $2$, we may split it into
two pieces $p=p_0\union p_1$ where $\dom(p_0)\of \gamma\cross
\alpha$ and $\dom(p_1)\of \gamma\cross
[\alpha,\gamma^\plusplus)$. It follows that $j(p)=j(p_0)\union
j(p_1)$, with the key point being that the domain of $j(p_1)$ is
disjoint from the domain of any element of
$j(\add(\gamma,\alpha))$.  In particular, $r^+$ is compatible
with $j(p_1)$. Since also $r^+\leq q\leq j(p_0)$, it follows that
$r^+$ is compatible with $j(p)$, as we claimed.

Now we iterate this idea to construct the
$M[j(G_\gamma)]$-generic filter $g^*\of
j(\add(\gamma,\gamma^\plusplus))$. Since the forcing
$j(\add(\gamma,\gamma^\plusplus))$ has size $j(\gamma^\plusplus)$
and is $j(\gamma^\plus)$-c.c., it has $j(\gamma^\plusplus)$ many
maximal antichains in $M[j(G_\gamma)]$. Since
$\card{j(\gamma^\plusplus)}^V=\gamma^\plusplus$, we may enumerate
these antichains in a sequence
$\<A_\beta\st\beta<\gamma^\plusplus>$ in $V[G_\gamma][g]$. We now
define a descending sequence $\<r_\beta\st
\beta<\gamma^\plusplus>$ of conditions in
$j(\add(\gamma,\gamma^\plusplus))$, each of which is compatible
with every element of $j\image g$. At successor stages, if
$r_\beta$ is defined, we employ the argument of the previous
paragraph to select a condition $r_{\beta+1}\leq r_\beta$ that
decides the antichain $A_\beta$ and remains compatible with
$j\image g$. At limit stages $\eta$, let
$r_\eta=\union\set{r_\beta\st \beta<\eta}$. Because
$M[j(G_\gamma)]$ is closed under $\gamma^\plus$ sequences in
$V[G_\gamma][g]$, we know that $r_\eta\in M[j(G_\gamma)]$, and it
clearly remains compatible with every element of $j\image g$.

Let $g^*$ be the filter generated by the conditions
$\set{r_\beta\st \beta<\gamma^\plusplus}$. Since these conditions
decide every maximal antichain, $g^*$ is
$M[j(G_\gamma)]$-generic. And since $j\image g\of g^*$, we may
lift the embedding to $j:V[G_\gamma][g]\to M[j(G_\gamma)][j(g)]$,
where $j(g)=g^*$, thereby witnessing the
$\gamma^\plus$-supercompactness of $\gamma$ in $V[G_\gamma][g]$,
as desired.

The careful reader will observe that we have actually proved that
if $\gamma$, a limit of partially supercompact cardinals, is
$\lambda$-supercompact in $V$ for some regular cardinal $\lambda$
above $\gamma$, then this is preserved by forcing over
$V[G_\gamma]$ with $\add(\gamma,\gamma^\plusplus)$. In
particular, the case $\gamma=\kappa$ shows that the
supercompactness of $\kappa$ in $V[G]$ is indestructible by
$\add(\kappa,\kappa^\plusplus)$, just as the theorem states, and
so the proof is complete.\QED

\Corollary. After the lottery preparation preserving
level-by-level agreement, the supercompactness of $\kappa$
becomes indestructible by\/ $\add(\theta,\theta^\plus)$ for any
regular cardinal\/ $\theta\geq\kappa$.\label{thetaplus}

\Proof: We use a similar argument for this Corollary. First we
claim for any $\gamma\leq\kappa$ that forcing with
$\add(\theta,\theta^\plus)$ over $V[G_\gamma]$, where
$\theta\geq\gamma$ is a regular cardinal below the next partially
supercompact cardinal above $\gamma$, preserves the level-by-level
agreement between strong compactness and supercompactness.
Suppose inductively that this holds below $\gamma$ and that
$g\of\add(\theta,\theta^\plus)$ is $V[G_\gamma]$-generic. It is
easy to see that the level-by-level agreement below or above
$\gamma$ is not affected, so suppose that $\gamma$ is
$\lambda$-strongly compact for some regular cardinal
$\lambda>\gamma$ in $V[G_\gamma][g]$; we will show that $\gamma$
is $\lambda$-supercompact there as well. By the Gap Forcing
Theorem \cite{GapForcing} we know that $\gamma$ is
$\lambda$-strongly compact and hence $\lambda$-supercompact in
$V$.

If $\lambda<\theta$, then since $\add(\theta,\theta^\plus)$ adds
no new subsets to $P_\gamma\lambda$, the
$\lambda$-supercompactness of $\gamma$ in $V[G_\gamma]$ is
trivially preserved to $V[G_\gamma][g]$. So we may assume
$\theta\leq\lambda$. In this case it follows that $\gamma$ is a
limit of partially supercompact cardinals, since otherwise the
forcing $\P_\gamma*\add(\theta,\theta^\plus)$ would be equivalent
to small forcing followed by $\ltgamma$-closed forcing adding a
subset to $\theta$, which by \cite{Dual} would destroy the
$\theta$-strong compactness of $\gamma$, contrary to our
assumption that $\gamma$ is $\lambda$-strongly compact in
$V[G_\gamma][g]$. So, let $j:V\to M$ be a
$\lambda$-supercompactness embedding such that $\gamma$ is not
$\lambda$-supercompact in $M$. Since the induction hypothesis
holds up to $j(\gamma)$ in $M$, the forcing
$\add(\theta,\theta^\plus)$ appears in the stage $\gamma$ lottery
of $j(\P_\gamma)$, and below a condition opting for this poset in
that lottery we may factor the iteration as
$\P_\gamma*\add(\theta,\theta^\plus)*\Ptail$, where $\Ptail$ is
$\leqlambda$-closed in $M[G_\gamma][g]$. And the diagonalization
argument of Corollary \ref{gchind} shows how to lift this
embedding to $j:V[G_\gamma]\to M[j(G_\gamma)]$, where
$j(G_\gamma)=G_\gamma*g*\Gtail$ for some generic filter
$\Gtail\of\Ptail$ constructed in $V[G_\gamma][g]$. It remains to
lift the embedding through the forcing
$\add(\theta,\theta^\plus)$.

If $\theta^\plus\leq\lambda$, this can be done with the usual
master condition argument, and we omit the details. The remaining
case, the hard case of $\theta=\lambda$, proceeds as in the hard
case of Corollary \ref{gchind}. Specifically, one first shows as
before that if $r\in j(\add(\theta,\theta^\plus))$ is compatible
with $j\image g$ and $A\of j(\add(\theta,\theta^\plus))$ is a
maximal antichain in $M[j(G_\gamma)]$, then there is a stronger
condition $r^+\leq r$ deciding $A$ and still compatible with
$j\image g$. For this, one uses the fact that
$j\image\theta^\plus$ is unbounded in $j(\theta^\plus)$ and
consequently $A$ is contained in $j(\add(\theta,\alpha))$ for
sufficiently large $\alpha<\theta^\plus$. By counting antichains
and iterating this argument, we once again construct a descending
sequence of conditions that eventually meet every maximal
antichain of $j(\add(\theta,\theta^\plus))$ in $M[j(G_\gamma)]$.
The filter $g^*$ generated by these conditions is therefore
$M[j(G_\gamma)]$-generic and extends $j\image g$, so we may lift
the embedding to $j:V[G_\gamma][g]\to M[j(G_\gamma)][j(g)]$,
where $j(g)=g^*$, thereby witnessing the
$\lambda$-supercompactness of $\gamma$ in $V[G_\gamma][g]$, as
desired.

Finally, we observe that we have actually proved that if $\gamma$
is a limit of partially supercompact cardinals and is
$\lambda$-supercompact in $V$, then forcing with
$\add(\theta,\theta^\plus)$ over $V[G_\gamma]$ for any
$\theta\leq\lambda$ preserves the $\lambda$-supercompactness of
$\gamma$. In particular, the case $\gamma=\kappa$ shows that the
full supercompactness of $\kappa$ is indestructible by
$\add(\theta,\theta^\plus)$ for any regular
$\theta\geq\kappa$.\QED

We regret that we have not been able to generalize these results
to $\add(\kappa,\kappa^\plusplusplus)$ or
$\add(\theta,\theta^\plusplus)$ for $\theta>\kappa$. Nevertheless,
because the results of this section show that after the lottery
preparation preserving level-by-level agreement the supercompact
cardinal $\kappa$ becomes indestructible by a great variety of
forcing notions, and we know of no specific $\ltkappa$-directed
closed forcing notion which does not preserve the
supercompactness of $\kappa$ over this model, we regard it as
currently the most natural candidate for a positive answer to
Question \ref{Open Question}, if any exists. The key questions
then become:

\Question. After the lottery preparation preserving
level-by-level agreement, is the supercompact cardinal $\kappa$
fully indestructible? If not, for which posets does it become
indestructible?

To conclude this section, we will show that the lottery
preparation preserving level-by-level agreement makes the
supercompact cardinal resurrectible after any $\ltkappa$-directed
closed forcing. Recall that a supercompact cardinal $\kappa$ is
said to be {\df resurrectible} if after any $\ltkappa$-directed
closed forcing $\Q$ there is further $\ltkappa$-distributive
forcing $\R$ such that $\Q*\R$ preserves the supercompactness of
$\kappa$. Thus, even if $\Q$ destroys the supercompactness of
$\kappa$, it is recovered by further forcing with $\R$.

\Corollary. After the lottery preparation preserving
level-by-level agreement, the supercompactness of\/ $\kappa$ is
resurrectible after any $\ltkappa$-directed closed forcing.

\Proof: The point here is that indestructibility by
$\coll(\kappa,\theta)$ implies resurrectibility by any
$\ltkappa$-directed closed forcing $\Q$ of size $\theta$,
assuming $\theta^{\ltkappa}=\theta$. This is true because such a
poset $\Q$ completely embeds into the collapse poset (this can be
seen by observing that $\Q\cross\coll(\kappa,\theta)$ is
$\ltkappa$-closed, collapses $\theta$ to $\kappa$ and has size
$\theta$, and noting that there is only one forcing notion with
these features). Thus, we may view $\coll(\kappa,\theta)$ as
$\Q*\R$, where $\R$ is $\coll(\kappa,\theta)$, and so even if $\Q$
happens to destroy the supercompactness of $\kappa$, further
forcing with $\R$ amounts altogether to forcing with
$\Q*\R=\coll(\kappa,\theta)$, which preserves the supercompactness
of $\kappa$.\QED

The previous argument actually provides a stronger kind of
resurrectibility than we claimed. Namely, define that $\kappa$ is
{\df $\theta$-resurrectible} if for any $\lttheta$-closed $\Q$
there is $\lttheta$-distributive forcing $\R$ such that $\Q*\R$
preserves the supercompactness of $\kappa$. Thus, for example,
$\kappa$ is resurrectible if and only if it is
$\kappa$-resurrectible. We now define that $\kappa$ is {\df
strongly resurrectible} if it is $\theta$-resurrectible for every
$\theta\geq\kappa$. Since any $\lttheta$-closed forcing $\Q$
embeds completely into $\coll(\theta,\card{\Q})$, and the
supercompactness of $\kappa$ is indestructible by such forcing,
we have actually proved:

\Corollary. After the lottery preparation preserving
level-by-level agreement, the supercompactness of\/ $\kappa$
becomes strongly resurrectible.\label{StronglyResurrectible}

If indeed all one desires is resurrectibility, then actually
there is no need as in Theorem \ref{Near} for restricting the
possibility of large cardinals above the supercompact cardinal in
question. Specifically, we claim the following:

\Theorem. If $\kappa$ is the least supercompact cardinal---whether
or not there are large cardinals above $\kappa$---there is a
forcing extension, preserving all supercompact cardinals, with a
level-by-level agreement between strong compactness and
supercompactness, in which the supercompactness of $\kappa$
becomes strongly resurrectible.

\Proof: We may assume, by forcing with the poset of
\cite{ApterShelah495} if necessary, that there is a level-by-level
agreement between strong compactness and supercompactness in $V$
and that the \GCH\ holds there (note that the
\cite{ApterShelah495} forcing preserves all supercompact
cardinals and by the Gap Forcing Theorem \cite{GapForcing} creates
no new ones). Because $\kappa$ is supercompact, it is a limit of
strong cardinals (see \cite[Lemma 2.1 and the subsequent
remark]{AptCum2}). Furthermore, since $\kappa$ is the least
supercompact cardinal in $V$, no cardinal below $\kappa$ is
supercompact up to a strong cardinal cardinal (lest it be fully
supercompact). Let $\P$ be the reverse Easton support
$\kappa$-iteration with nontrivial forcing only at stages
$\gamma<\kappa$ that are inaccessible limits of strong cardinals.
At such a stage $\gamma$, the forcing is the lottery sum of all
$\coll(\theta_1,\theta_2)$, where
$\gamma\leq\theta_1\leq\theta_2$, each $\theta_i$ is regular, and
$\theta_2$ is less than the next strong cardinal above $\gamma$
(plus trivial forcing). Suppose that $G\of\P$ is $V$-generic.

The usual lifting arguments (see e.g. \cite{LotteryPreparation})
establish that $\kappa$ remains supercompact in $V[G]$ and
furthermore that the supercompactness of $\kappa$ becomes
indestructible there by further forcing with
$\coll(\theta_1,\theta_2)$ whenever
$\kappa\leq\theta_1\leq\theta_2$ and each $\theta_i$ is regular.
(Note: the possibility of strong cardinals above $\kappa$ is
irrelevant here, since for any $\theta$ one can use an embedding
$j:V\to M$ for which $\kappa$ is not $\theta$-supercompact in
$M$; consequently, the next strong cardinal in $M$ above $\kappa$
is above $\theta$, which is all that is needed in the lifting
argument). It follows as in Corollary \ref{StronglyResurrectible}
that the supercompactness of $\kappa$ is strongly resurrectible in
$V[G]$.

We now argue that the model $V[G]$ retains the level-by-level
agreement between strong compactness and supercompactness. Since
the forcing $\P$ has size $\kappa$, we know that the
level-by-level agreement for cardinals above $\kappa$ holds
easily by \cite{LevySolovay67}. It remains only to consider
cardinals $\gamma$ below $\kappa$. Accordingly, suppose that
$\gamma$ is $\theta$-strongly compact in $V[G]$ for some regular
cardinal $\theta>\gamma$; we aim to show it is also
$\theta$-supercompact there. By the Gap Forcing Theorem
\cite{GapForcing}, we know that $\gamma$ is $\theta$-strongly
compact in $V$, and hence also $\theta$-supercompact there. Fix a
$\theta$-supercompactness embedding $j:V\to M$ for which $\gamma$
is not $\theta$-supercompact in $M$. It follows that the next
strong cardinal in $M$ above $\gamma$ is above $\theta$.

We first treat the case in which $\gamma$ is a limit of strong
cardinals, that is, when $\gamma$ is a stage of forcing in $\P$.
In the stage $\gamma$ lottery, the generic $G$ selected some
winning poset $\Q$, and below a condition deciding this we may
factor the forcing $\P$ as $\P_\gamma*\Q*\Ptail$ where $\Q$ is
either trivial forcing or $\coll(\theta_1,\theta_2)$ for some
$\gamma\leq\theta_1\leq\theta_2$. Since $\Ptail$ is closed beyond
the next strong cardinal above $\gamma$, it does not affect the
$\theta$-supercompactness of $\gamma$, and so it suffices for us
to show that $\gamma$ is $\theta$-supercompact in
$V[G_\gamma][g]$, where $g\of\Q$ is $V[G_\gamma]$-generic. If
$\Q$ is trivial, or if $\theta_2\leq\theta$, then the usual
lifting arguments allow us to lift the embedding $j:V\to M$ to
$j:V[G_\gamma][g]\to M[j(G_\gamma)][j(g)]$, thereby witnessing
the $\theta$-supercompactness of $\gamma$ in $V[G_\gamma][g]$. So
we may assume $\theta<\theta_2$. Since $\theta$ is a cardinal in
$V[G]$, it follows that $\theta\leq\theta_1$. If
$\theta<\theta_1$, then the forcing $\Q$ adds no new subsets to
$P_\gamma\theta$, and so the $\theta$-supercompactness of $\gamma$
in $V[G_\gamma]$ (from the case of trivial $\Q$ above) is
preserved to $V[G_\gamma][g]$. So we have reduced to the case
that $\theta=\theta_1$. Since $\gamma$ is $\theta$-strongly
compact in $V[G]$, it is also $\theta_2$-strongly compact there,
and hence $\theta_2$-supercompact in $V$. One may therefore
employ the usual arguments to lift a $\theta_2$-supercompactness
embedding from $V$ to $V[G_\gamma][g]$, as desired.

We now treat the case that $\gamma$ is not a limit of strong
cardinals. Let $\delta<\gamma$ be the supremum of the strong
cardinals below $\gamma$. If $\delta$ is not inaccessible, then
the forcing $\P$ factors as $\P_\delta*\Ptail$, where $\P_\delta$
is small relative to $\gamma$ and $\Ptail$ is closed beyond the
next strong cardinal above $\gamma$. Such forcing must preserve
the $\theta$-supercompactness of $\gamma$. We may therefore
assume alternatively that $\delta$ is inaccessible, and hence a
stage of forcing. The generic $G$ selected a winning poset $\Q$
in the stage $\delta$ lottery, and below a condition deciding
this we may factor $\P$ as $\P_\delta*\Q*\Ptail$. The forcing
$\Ptail$ is closed beyond $\theta$, and does not affect the
$\theta$-supercompactness of $\gamma$. Thus, it suffices for us
to see that $\gamma$ is $\theta$-supercompact in
$V[G_\delta][g]$, where $g\of\Q$ is $V[G_\delta]$-generic. The
forcing $\Q$ is either trivial or $\coll(\theta_1,\theta_2)$ for
some $\delta\leq\theta_1\leq\theta_2$. If $\Q$ is trivial or
$\theta_2<\gamma$, then the forcing is small relative to
$\gamma$, and the result is immediate. So assume
$\gamma\leq\theta_2$. Since $\gamma$ is a cardinal, it must be
that $\gamma\leq\theta_1$. If in addition $\theta<\theta_1$, then
$\gamma<\theta_1$ and we may ignore the forcing $\Q$, since it
does not destroy the $\theta$-supercompactness of $\gamma$ in
$V[G_\delta]$, a small forcing extension. What remains is the case
$\gamma\leq\theta_1\leq\theta$. Here, as in Corollary
\ref{thetaplus}, we make a key use of the main result of
\cite{Dual}: the forcing $\P_\delta*\coll(\theta_1,\theta_2)$ is
small forcing followed by $\ltgamma$-closed forcing that adds a
new subset to $\theta_1$. By \cite{Dual}, such forcing
necessarily destroys the $\theta_1$-strong compactness of
$\gamma$, contrary to our assumption that $\gamma$ is
$\theta$-strongly compact in $V[G]$, and hence in
$V[G_\delta][g]$. So the proof is complete.\QED

Because of this argument, the limitation identified above the
supercompact cardinal in Theorem \ref{No} does not engage for
resurrectibility.

\Section Indestructibility with near level-by-level agreement
\label{NearSec2}

In this section, we will prove that full indestructibility is
compatible with a level-by-level agreement between strong
compactness and supercompactness almost everywhere. To make this
notion of almost everywhere precise, let us define that a set
$A\of\kappa$ is {\df large} with respect to supercompactness if
for every $\theta\geq\kappa$ there is a $\theta$-supercompactness
embedding $j:V\to M$ with $\kappa\in j(A)$. That is, the set is
large in virtue of having measure one with respect to all these
induced normal measures.\footnote{We warn the reader that the
collection of large subsets of a supercompact cardinal do not
form a filter. Specifically, if $\cal M_\theta$ is the collection
of normal measures induced by $\theta$-supercompactness embeddings
with critical point $\kappa$, then the collection of large sets is
precisely $\intersect_\theta(\union\cal M_\theta)$. Since $\cal
M_\theta\of\cal M_\lambda$ whenever $\lambda\leq\theta$, the
collections $\cal M_\theta$ are eventually equal, and the
collection of large sets is simply the union of the measures in
this stabilized eventual value for $\cal M_\theta$. Since the
proof of Example \ref{Robust} shows that each $\cal M_\theta$ has
many measures, this means that the collection of large sets is
the union of a great number of normal measures, and therefore it
cannot be a filter.}

To illustrate this notion, define that a measurable cardinal
$\gamma$ is {\df robust} if whenever $\gamma$ is
$\theta$-supercompact, then there is a $\theta$-supercompactness
embedding $j:V\to M$ with $\cp(j)=\gamma$ such that $\gamma$ is
$\theta$-supercompact in $M$. Equivalently, $\gamma$ is robust if
it has nontrivial Mitchell rank in every degree of
supercompactness that it exhibits. Conversely, we could define
that $\gamma$ is {\df precarious} when it is $\eta$-supercompact
with trivial Mitchell rank for some $\eta$; these cardinals
necessarily lose of bit of their supercompactness in their most
supercompact ultrapowers. Thus, by definition, the robust
cardinals are exactly the non-precarious measurable cardinals.

\Example. There are many robust cardinals below any supercompact
cardinal; indeed, the collection of robust cardinals is large
with respect to supercompactness.\footnote{The set of precarious
cardinals is also large with respect to supercompactness, since
by a suitable choice of supercompactness embedding---using a
measure of Mitchell rank $1$---one can arrange that $\kappa$ has
trivial Mitchell rank in $M$ for its largest degree of
supercompactness.}\label{Robust}

\Proof: The basic point is that if a cardinal $\kappa$ is
$2^{\theta^\ltkappa}$-supercompact, then it has very high
Mitchell rank in $\theta$-supercompactness. To see this, suppose
$j:V\to M$ is a $2^{\theta^\ltkappa}$-supercompactness embedding.
Since all the $\theta$-supercompactness measures from $V$ are in
$M$, it follows that $\kappa$ is $\theta$-supercompact in $M$.
Thus, for the induced $\theta$-supercompactness factor embedding
$j_0:V\to M_0$, the cardinal $\kappa$ is $\theta$-supercompact in
$M_0$. So the Mitchell rank is at least $1$ in $V$, and so it is
at least $1$ in $M$ and therefore also in $M_0$; so it is at
least $2$ in $V$, and so at least $2$ in $M$ and therefore at
least $2$ in $M_0$, and so on cycling around up to $\theta$ and
beyond.

To prove the claim, now, let $\theta$ be any strong limit
cardinal above $\kappa$ and let $j:V\to M$ be a
$\theta$-supercompactness embedding by a measure with trivial
Mitchell rank (or just make $j(\kappa)$ as small as possible). It
follows that $\kappa$ is not $\theta$-supercompact in $M$; but by
the closure of $M$ it is $\lttheta$-supercompact there.
Therefore, since $\theta$ is a strong limit cardinal, the basic
point in the previous paragraph shows that it exhibits nontrivial
Mitchell rank in $M$ for every degree of supercompactness below
$\theta$. So it is robust in $M$. And since this is true for
arbitrarily large $\theta$, it follows that the set of robust
cardinals below $\kappa$ is large with respect to
supercompactness.\QED

While the argument in the example verifies the robustness of a
cardinal by having it exhibit a limit degree of supercompactness,
this is not at all the only way to be robust, for a cardinal
$\kappa$ could be only $\kappa^\plus$-supercompact, for example,
and still be robust by simply having nontrivial Mitchell degree
for $\kappa^\plus$-supercompactness. The collection of cardinals
with a largest degree of supercompactness and nontrivial Mitchell
rank for that degree of supercompactness is also large with
respect to supercompactness; this can be seen by using
$\theta$-supercompactness embeddings $j:V\to M$ for $\theta$ a
successor cardinal arising from a $\theta$-supercompactness
measure with Mitchell rank at least $2$.

Let us now turn to the main theorem of this section:

\Theorem. Suppose that $\kappa$ is supercompact and no cardinal
is supercompact up to a larger cardinal $\lambda$ which is
$2^\lambda$-supercompact. Then there is a forcing extension in
which $\kappa$ becomes indestructibly supercompact and
level-by-level agreement holds between strong compactness and
supercompactness on a set that is large with respect to
supercompactness.

\Proof: As in Theorem \ref{Near}, we may assume without loss of
generality, by forcing if necessary, that the \GCH\ holds and
further, by forcing with the notion in \cite{ApterShelah495} if
necessary, that in $V$ there is already a level-by-level agreement
between strong compactness and supercompactness. Thus, we also
have that no cardinal is supercompact up to a partially
supercompact cardinal. Let $\P$ be the reverse Easton support
$\kappa$-iteration which adds a Cohen real and then has
nontrivial forcing only at cardinals $\gamma$ that are
inaccessible limits of partially supercompact cardinals. At such
a stage $\gamma$, the stage $\gamma$ forcing $\Q_\gamma$ is the
lottery sum of all $\ltgamma$-directed closed $\Q$ of size less
than $\theta_\gamma$, the least cardinal such that $\gamma$ is
not $\theta_\gamma$-supercompact in $V$. Suppose that $G\of\P$ is
$V$-generic, and consider the model $V[G]$.

First, we claim that $\kappa$ is indestructibly supercompact in
$V[G]$. It suffices to argue that if $\Q \in V[G]$ is
$\ltkappa$-directed closed in $V[G]$ and $H\of\Q$ is
$V[G]$-generic, then $\kappa$ is supercompact in $V[G][H]$. Fix
any regular cardinal $\theta > \card{\Q}$ and a
$\theta$-supercompactness embedding $j:V\to M$ so that $\kappa$
is not $\theta$-supercompact in $M$. Since $\kappa$ will
necessarily be $\lttheta$-supercompact in $M$, the forcing $\Q$
will appear in the stage $\kappa$ lottery of $j(\P)$. Below a
condition opting for $\Q$ in this lottery, the forcing $j(\P)$
factors as $\P*\Q*\Ptail$, where $\Ptail$ is the remainder of the
iteration. Note that since no cardinal is supercompact beyond a
partially supercompact cardinal, no cardinal above $\kappa$ is
partially supercompact. Thus, the next nontrivial stage of
forcing in $j(\P)$ is beyond $\theta$, and so $\Ptail$ is
$\leqtheta$-closed in $M[G][H]$. Therefore, by the usual
diagonalization techniques, we may construct in $V[G][H]$ an
$M[G][H]$-generic filter $\Gtail\of\Ptail$ and lift the embedding
to $j:V[G]\to M[j(G)]$ with $j(G)=G*H*\Gtail$. After this, we
find a master condition below $j\image H$ in $M[j(G)]$ and again
by diagonalization construct an $M[j(G)]$-generic $j(H)\of j(\Q)$
in $V[G][H]$. This allows us to lift the embedding fully to
$j:V[G][H]\to M[j(G)][j(H)]$, which witnesses the
$\theta$-supercompactness of $\kappa$ in $V[G][H]$, as desired.

Second, we claim that in $V[G]$ there is a level-by-level
agreement between strong compactness and supercompactness at the
cardinals $\gamma$ that are robust in $V$. To see this, on the
one hand, the main results of \cite{Dual} show that since the
forcing $\P$ is mild and admits a very low gap (see
\cite{GapForcing} for the relevant definitions), it does not
increase the degree of strong compactness or supercompactness of
any cardinal. On the other hand, we will now argue that it fully
preserves all regular degrees of supercompactness of every
partially supercompact cardinal $\gamma$ that is robust in $V$.
Suppose that $\gamma$ is robust and $\theta$-supercompact in $V$
for some regular cardinal $\theta$ above $\gamma$. A simple
reflection argument shows that $\gamma$ must be a limit of
partially supercompact cardinals, and therefore is a nontrivial
stage of forcing. The generic filter $G$ opted for some
particular forcing $\Q$ in the stage $\gamma$ lottery. By
increasing $\theta$ if necessary, we may assume that
$\card{\Q}\leq\theta$. Since no cardinal is supercompact beyond a
partially supercompact cardinal, as in Lemma \ref{SubA}, the next
nontrivial stage of forcing beyond $\gamma$ is well beyond
$\theta$. It therefore suffices to argue that $\gamma$ is
$\theta$-supercompact in $V[G_{\gamma+1}]$. Fix any
$\theta$-supercompactness embedding $j:V\to M$ with critical
point $\gamma$ such that $\gamma$ is $\theta$-supercompact in $M$
(this is where we use the robustness of $\gamma$). It follows
that $\Q$ appears in the stage $\gamma$ lottery of
$j(\P_\gamma)$, and we may lift the embedding to
$j:V[G_{\gamma+1}]\to M[j(G_{\gamma+1})]$ just as we did in the
previous paragraph with $\kappa$. Thus, $\gamma$ remains
$\theta$-supercompact in $V[G_{\gamma+1}]$ and hence in $V[G]$.
So we retain the level-by-level agreement between strong
compactness and supercompactness for the partially supercompact
cardinals $\gamma$ that are robust in $V$.

Now let us argue that the collection $A$ of partially
supercompact cardinals that are robust in $V$ remains large with
respect to supercompactness in $V[G]$. The point is that if
$j:V\to M$ is a $\theta$-supercompactness embedding with critical
point $\kappa$, then below a condition opting for trivial forcing
in the lottery at stage $\kappa$, we may lift the embedding to
$j:V[G]\to M[j(G)]$. In particular, if $\kappa\in j(A)$ for the
original embedding, then this remains true for the lifted
embedding. Indeed, this argument shows that our forcing preserves
every set that is large with respect to supercompactness in the
ground model.\QED

We would like to point out that in fact in $V[G]$ we have much
more level-by-level agreement than just at the robust cardinals.
In particular, for any cardinal $\gamma$ that is a limit of
partially supercompact cardinals, if it happens that the generic
filter chooses a poset $\Q$ with $\card{\Q}^\plus<\theta_\gamma$,
and this happens on a large set since we can arrange it at stage
$\kappa$ in $j(\P)$, then one may apply the lifting argument of
the proof using an embedding witnessing enough supercompactness
of $\gamma$.

The problematic cardinals are exactly the precarious cardinals
$\gamma$ for which the forcing opts for $\Q$ at stage $\gamma$ of
the largest possible size. This case is difficult because $\Q$
will be too large to appear in the stage $\gamma$ lottery of
$j(\P)$, and so one may not use any ordinary lifting argument to
preserve the supercompactness of $\gamma$. Despite our many
attempts, there seems to be no easy solution to this annoying
problem: restricting the posets that appear in the stage $\gamma$
lottery on the $V$-side simply causes a corresponding restriction
on the $j$-side, and it seems that there are always these
borderline posets that are allowed in the $\P$ lottery but not in
the $j(\P)$ lottery.

Even for these problematic cardinals $\gamma$, though, it will
often happen that nevertheless the level-by-level agreement
between strong compactness and supercompactness is preserved,
because the forcing $\Q$ opted for in the stage $\gamma$ lottery
is, say, equivalent to smaller forcing, or does not add sets below
$\card{\Q}$ but destroys the $\card{\Q}$-supercompactness of
$\gamma$ (so that both strong compactness and supercompactness
drop evenly). Nevertheless, we know of no way to ensure the
level-by-level agreement at all cardinals below $\kappa$, so
Question \ref{Open Question} remains open.

The theorems of this section, as well as the previous section,
suggest that a positive answer to Question \ref{Open Question} is
possible. Further, since the cardinals below $\kappa$ at which
the level-by-level agreement holds exhibit themselves some degree
of indestructibility, the results suggest the intriguing
possibility that one could have a supercompact cardinal and a
level-by-level agreement between strong compactness and
supercompactness in the presence of {\df universal
indestructibility}: every partially supercompact cardinal
$\gamma$ is fully indestructible by $\ltgamma$-directed closed
forcing.

\Question. Is universal indestructibility consistent with
level-by-level agreement if there is a supercompact cardinal?

Universal indestructibility, like the level-by-level agreement
between strong compactness and supercompactness in Theorem
\ref{No}, is by itself incompatible with large cardinals above a
supercompact cardinal. Specifically, in
\cite{UniversalIndestructibility} we show that when universal
indestructibility holds, then no cardinal is supercompact beyond
a measurable cardinal. The question is, does this affinity
indicate a compatibility of the two notions?

\bibliographystyle{alpha}
\bibliography{MathBiblio}

\end{document}